\documentclass[a4paper,12pt]{article}
\usepackage[T2A]{fontenc}
\usepackage[utf8]{inputenc}
\usepackage[english, russian]{babel}
\usepackage{amsmath,amsthm}
\usepackage{amssymb}
\usepackage{amscd}
\usepackage{amsfonts}

\usepackage[matrix,arrow,curve]{xy}

\theoremstyle{plain}
\newtheorem{theo}{Теорема}[section]

\newtheorem{lem}[theo]{Лемма}

\newtheorem{defi}[theo]{Определение}
\newtheorem{theod}[theo]{Теорема}

\renewcommand{\lim}{\varprojlim}

\renewcommand{\qed}{{\hfill$\square$}}

\begin{document}
\selectlanguage{russian}

\title{\normalsize{Циклические накрытия, которые не являются стабильно рациональными\\ $\,$\\ Cyclic covers that are not stably rational}}

\author{\normalsize{Жан-Луи Кольё-Телэн}\thanks{
Jean-Louis Colliot-Th\'el\`ene,
C.N.R.S., Universit\'e Paris Sud, Math\'ematiques, B\^atiment 425, 91405 Orsay Cedex, France; 
Chaire Lam\'e 2015,  Universit\'e d'\'Etat de Saint P\'etersbourg, Исследовательская Лаборатория имени П. Л. Чебышёва, Saint-P\'etersbourg, Russie;  jlct@math.u-psud.fr }
$\,$ \normalsize{и}  
$\,$  \normalsize{Елена Пирютко} \thanks{
Alena Pirutka,
C.N.R.S., \'Ecole Polytechnique, CMLS,  91128 Palaiseau, France; 
alena.pirutka@polytechnique.edu
}
}

\date{}

\maketitle

\vspace{-0.8 cm}

 \begin{abstract}
На основе методов, разработанных Колларом, Вуазан, авторами, Тотаро, 
мы доказываем, что циклическое накрытие $\mathbb P_{\mathbb C}^n, n\geq 3$  простой степени $p$,  разветвлённое над очень общей гиперповерхностью 
 $f(x_0,\ldots , x_n)=0$ степени $mp$ не является стабильно рациональным при условии  
   $n+1\leq mp$.  В размерности $3$ получаем двойные накрытия $\mathbb P^3_{\mathbb C}$, разветвлённые над очень общей гиперповерхностью степени $4$ (Вуазан),  а также   двойные накрытия $\mathbb P^3_{\mathbb C}$, разветвлённые над очень общей гиперповерхностью степени $6$ (Бовиль). Мы также получаем   двойные накрытия $\mathbb P^4_{\mathbb C}$, разветвлённые над очень общей гиперповерхностью степени $6$.  Метод статьи позволяет получить примеры над числовыми полями.\\
{\it Ключевые слова }:  стабильная рациональность, группа Чжоу нуль-циклов, циклические накрытия.

\begin{center}{\bf R\'esum\'e}\\
\end{center}

Using the methods developed by Koll\'ar, Voisin, ourselves, Totaro, we prove that a cyclic cover of   $\mathbb P_{\mathbb C}^n, n\geq 3$  of prime degree $p$, ramified along a very general hypersurface of degree $mp$ is not stably rational if  $n+1\leq mp$. 
In small dimensions, we recover double covers of  $\mathbb P^3_{\mathbb C},$  ramified along a quartic (Voisin), and double covers of $\mathbb P^3_{\mathbb C}$ ramified along a sextic (Beauville), and we also find double covers of $\mathbb P^4_{\mathbb C}$ ramified along a sextic. This method also allows one to produce examples over a number field. \\
   {\it Keywords}: stable rationality, Chow group of zero-cycles, cyclic covers.\\
  \end{abstract}
  
 \noindent Индекс УДК : 512.752


\section{Введение}

 Проективное многообразие  $X$  над полем $k$ называется стабильно рациональным,  если для некоторого $n$ многообразие $X\times\mathbb P^n_k$ является рациональным. Существуют стабильно рациональные нерациональные многообразия \cite{BCTSSD}.  В  работе \cite{V13} Клэр Вуазан вводит метод для доказательства, что многообразие  $X$ не является стабильно рациональным, основанный на целом разложении диагонали в группе Чжоу $CH^{dim\, X}(X\times X)$ и специализации. Этот метод позволяет доказать, что двойное накрытие $\mathbb P^3_{\mathbb C}$, разветвлённое над очень общей поверхностью степени $4$, не является стабильно рациональным.   В работе \cite{CTP15} мы рассматриваем свойство $CH_0$- универсальной тривиальности, эквивалентное целому разложению диагонали для гладких проективных многообразий, и которое делает метод специализации более гибким, в частности, для случая специализации над кольцом дискретного нормирования с полем вычетов положительной характеристики. 
  Mы доказываем, что при очень общем выборе коэффициентов,  гладкая комплексная квартика размерности $3$ не является стабильно рациональным многообразием.\\
 
  {\defi{ \it Пусть $f:X\to Y$ проективный морфизм  
 многобразий над полем $k$. Морфизм $f$ называется $CH_0$- универсально тривиальным, если для любого расширения полей $L/k$ отображение $f_*:CH_0(X_L)\to CH_0(Y_L)$  является изоморфизмом.  Если $Y=Spec\,k$ и $f$ -- структурный морфизм, то  многообразие $X$  называется $CH_0$- универсально тривиальным. }}
 
 В частности, стабильно рациональное многообразие является $CH_0$- универсально тривиальным.\\
 
  В работах А.  Бовиля \cite{B4, B6} рассматриваются случаи двойных накрытий $\mathbb P^3_{\mathbb C}$, разветвлёных над очень общей поверхностью степени $6$, а также двойных накрытий $\mathbb P^4_{\mathbb C}$ и  $\mathbb P^5_{\mathbb C}$, разветвлёных над очень общей гиперповерхностью степени $4$.
В работе \cite{HKT15},  Э. Креш, Б. Хассетт и Ю. Чинкель рассматривают случай некоторых расслоений на коники.

   Б. Тотаро \cite{To15}  доказал, что очень общая поверность степени $d$ в $\mathbb P^{n+1}_{\mathbb C}$ не является стабильно рациональной при условиях   
    $d\geq 2\lceil (n+2)/3 \rceil$  и  $n\geq 3$; в доказательстве  используются   результаты Коллара \cite{Ko96, KSC04} o двойныx накрытияx в характеристике $2$ и результат о специализации $CH_0$-универсальной тривиальности \cite{CTP15} 1.14 над кольцом дискретного нормирования с полем частных характеристики ноль и с полем вычетов положительной характеристики. Как замечает Тотаро \cite{To15}, методы, описанные выше, также возможно применить для более общих накрытий :    
    в этой работе мы продолжаем 
     изучение циклический накрытий  в положительной характеристике и  доказываем следующий результат  (теорема \ref{revnsr}):
 
{\theod{ \it Пусть $p$ -- простое число.  Пусть $X$ -- циклическое накрытие $\mathbb P_{\mathbb C}^n, n\geq 3$ степени $p$,   разветвлённое над очень общей гиперповерхностью $f(x_0,\ldots , x_n)=0$ степени $mp$.   Предположим,  что $m(p-1)  <n+1\leq mp$.   Тогда $X$  -- многообразие Фано, которое не является стабильно рациональным  многообразием. \\}}
 
 Как и в работе  \cite{To15}, мы получаем примеры над числовыми полями.
 
 Заметим, что  при $n=3, m=p=2$  получаем очень общие двойные 
  накрытия
   $\mathbb P^3_{\mathbb C}$, разветвлённые над квартикой (более общие результаты получены в работе \cite{V13}),   при $n=3, m=3, p=2$  мы получаем другое доказательство результатов \cite{B6}. 
 
 При условии $mp>n+1$, Коллар доказал, что рассматриваемые накрытия не являются линейчатыми \cite{Ko95}. Однако, это не даёт результатов о стабильной рациональности, так как существуют стабильно рациональные многообразия размерности $3$, которые не являются рациональными \cite{BCTSSD}.

 \section{$CH_0$-унивесальная тривиальность сингулярных многообразий}

\lem{Пусть $k$ -- алгебраически замкнутое поле и пусть $X$ -- целое проективное многообразие над $k$. Пусть $U\subset X$ открытое подмногообразие. Тогда для любой точки $z\in X(k)$ существует цикл $\xi\in Z_0(U)$ рационально эквивалентный  $z$ в $CH_0(X)$.}
\proof{Если $X=C$ -- целая кривая с нормализацией $D$, то утверждение следует из того, что группа Пикара  полулокальных колец $D$ тривиальна. В общем случае достаточно заметить, что существует целая кривая $C$, такая, что $z\in C$ и $C\cap U\neq \emptyset$. \qed\\}

\lem{Пусть $k$ -- алгебраически замкнутое поле и пусть $X$ -- целое проективное   $k$-рациональное многообразие. Если $X$ является гладким за исключением конечного числа точек, то $X$ -- универсально $CH_0$-тривиальное многообразие.}
\proof{Пусть $\emptyset\neq U\subset X$ открытое подмногообразие, изоморфное открытому подмногообразию $\mathbb P^n_k$. Пусть  $F/k$  -- некоторое расширение полей.  Любая гладкая точка $z\in X_F(F)$ рациональнo эквивалентна в $X_F$ нуль-циклу из $Z_0(U_F)$. Из предыдущей леммы получаем, что это остаётся верным для любой $F$-точки $X$. Аналогично рассуждениям \cite{CTP15} 1.5 получаем, что каждый цикл в $Z_0(X_F)$ рационально эквивалентен циклу $Nx$, для некоторого $N$ и (фиксированного)
$ x \in U(k) \subset U(F) \subset X(F)$.
 
   \qed\\}

 \lem\label{singord}{Пусть $k$ -- алгебраически замкнутое поле и пусть $X$ -- связное проективное многообразие над $k$. Если каждая приведённая компонента $X$ является $k$-рациональным многообразием с изолированными сингулярными точками, то $X$ -- универсально $CH_0$-тривиальное многообразие.}
 \proof{Достаточно применить предыдущую лемму и \cite{CTP15} 1.3. \qed\\}

В следующем параграфе мы применим лемму \ref{singord} для исключительных дивизоров разрешения особенностей.  
Приведём также более общее утверждение для  объединения  $CH_0$-универсально тривиальных многообразий. Далее в этой статье нам понадобится только  лемма \ref{singord}.
 
  \lem\label{union}{Пусть $X$ -- проективное 
  приведённое 
    геометрически  связное многообразие над полем $k$ и пусть $X=\bigcup_{i=1}^{N} X_i$ разложение $X$ на неприводимые компоненты. Предположим, что
\begin{itemize}
\item[(i)] каждое из многообразий $X_i$ геометрически неприводимо
и является $CH_0$-универсально тривиальным;
\item[(ii)]  каждое из пересечений $X_i\cap X_j$ либо пусто, либо содержит $0$-цикл $z_{ij}$ степени $1$.
\end{itemize}
Тогда многообразие $X$ является $CH_0$-универсально тривиальным многообразием.}
\proof{Пусть $L/k$ расширение полей и пусть $z\in CH_0(X_L)$  класс цикла степени $0$.  Так как  $X$ -- геометрически  связное, то в дуальном графе геометрических компонент существует полный цикл : существует последовательность индексов $i_1,\ldots i_m$, $1\leq i_j\leq N$ (где $m$  может быть больше, чем $N$), такая, что $\{i_1,\ldots, i_m\}=\{1,\ldots, N\}$ и $X_{i_j, L}\cap X_{i_{j+1}, L}$ непусто для всех $1\leq j\leq m$. 

Можно разложить  $z=\sum z_{i_j}$, где $z_{i_j}\in CH_0(X_{{i_j}L})$ степени $d_j$, $\sum d_j=0$ (с произвольным выбором на пересечениях, также некоторые $z_{i_j}$ могут быть равны нулю). Тогда $z_{i_1}=d_1z_{i_{1}i_{2}L}$ в $ CH_0(X_{{i_1}L})$, откуда $z_{i_1}+z_{i_2}=(d_1+d_2)z_{i_{2}i_{3}L}$  в  $ CH_0(X_{{i_1}L}\cup X_{{i_2}L})$  и т.д., получаем   $z=\sum z_{i_j}=(\sum d_i)z_{i_{m-1},i_{ m}L}=0$ в $CH_0(X_{L})$. \qed\\}

{\it Замечание.}  Условие $(i)$ выполняется, если существует разрешение особенностей  $\pi_i:\tilde X_i\to X_i$, такое, что  $\tilde X_i$ является $CH_0$-универсально тривиальным многообразием и все (схематические) слои $\pi_i$  являются $CH_0$-универсально тривиальными многообразиями (см. \cite{CTP15}, Prop. 1.8.)

\vspace{1cm}

\section{Циклические накрытия и особенности}

Напомним вкратце некоторые свойства циклических накрытий \cite{Ko96}.V, \cite{KSC04}.

\noindent Пусть $p$ -- простое число и $f(x_0, \ldots, x_n)$ -- однородный многочлен степени $mp$ над полем $k$.
{\it Циклическое накрытие } $\mathbb P^n_k$, разветвлённое над $f(x_0, \ldots, x_n)=0$, определяется как подмногообразие $\mathbb P(m,1,1,\ldots ,1)$ заданное условием $$y^p-f(x_0, \ldots , x_n)=0.$$

\noindent  Если  $k$ -- поле характеристики $p$, то такое циклическое накрытие почти всегда
сингулярно, особенности соответствуют критическим точкам $f$.\\

{\defi{ {\it Критической точкой} многочлена $g(x_1,\ldots, x_n)$ над полем $k$ называется  точка $P$, такая, что $\partial g/\partial x_i(P)=0\, \forall i$. Критическая точка  $P$ многочлена $g$ называется {\it невырожденной}, если определитель $|\frac {\partial^2g}{\partial x_i\partial x_j}(P)|$ не равен нулю. Критической точкой однородного многочлена $f(x_0, \ldots, x_n)$ называется критическая точка одного из многочленов $f(x_0,\ldots, x_{i-1}, 1, x_{i+1}, \ldots , x_n)$.  \\ }}

Если $k$ -- поле характеристики $2$ и $n$ нечётно, то все критические точки многочлена $g\in k[x_1,\ldots, x_n]$ являются вырожденными.\\ 

{\defi{ Пусть $k$ -- поле характеристики $2$ и $n$ нечётно.  Критическая точка  $P$ многочлена $g(x_1,\ldots , x_n)$ называется {\it почти невырожденной}, если $$\mathrm{length}\mathcal O_{\mathbb A^n, P}/(\partial g/\partial x_1(P), \ldots , \partial g/\partial x_n(P))=2.$$    }}

Для изучения стабильной рациональности нам понадобятся результаты о разрешении особенностей циклических накрытий (см. также \cite{Ko95}).

\lem\label{mortr}{Пусть $k$ -- алгебраически замкнутое поле характеристики $2$ и пусть $$X: \; y^2=f(x_1,\ldots, x_n)$$
аффинное  двойное накрытие, сингулярное в точке $P=(y, x_1,\ldots, x_n)=(0, 0,\ldots, 0)$, где  $n\geq 2$ чётно и $(0,\ldots, 0)$ -- невырожденная критическая точка многочлена $f$.  Пусть $\tilde X\to X$  -- раздутие точки $P$. 
Тогда :
\begin{itemize}
\item[(i)]   в окрестности точки $P$ накрытие $X$ задаётся условием 
$$y^2=x_1x_2+x_3x_4+\ldots + x_{n-1}x_n+g(x_1,\ldots, x_n),$$ 
где  каждый одночлен в записи $g(x_1,\ldots,  x_n)$ имеет степень не менее трёх.
\item[(ii)]  многообразие $\tilde X$ является гладким в окрестности  исключительного дивизора $E$ и  многообразие $E$ является универсально $CH_0$-тривиальным.
\end{itemize}  }
\proof{Свойство $(i)$ рассматривается в упражнении V.5.6.6 книги \cite{Ko96} (см. также доказательство теоремы \ref{revcp} далее). 

 Докажем $(ii)$. Достаточно рассмотреть следующие карты раздутия:
\begin{enumerate}
\item  $x_i=y z_i$, $1\leq i \leq n$ и $\tilde X$ задаётся условием
$$1= z_1 z_2+ z_3 z_4+\ldots +  z_{n-1} z_{n}+\frac{1}{y^2}g(yz_1,\ldots, y z_n)$$ в аффинных координатах $y, z_1, \ldots z_n$. 

Заметим, что многочлен $\frac{1}{y^2}g(y z_1,\ldots, y z_n)$ делится на $y$. 
Исключительный дивизор $E$ раздутия $\tilde X\to X$ задаётся условием $y=0$. Получаем уравнение $E$ в этой карте :
$$ z_1 z_2+ z_3 z_4+\ldots + z_{n-1} z_{n}=1,$$ задающее гладкую квадрику. 
Многообразие $\tilde X$ является гладким в каждой точке $E$ (это следует из уравнения $\tilde X$).

\item $y=wx_1$,  $x_i=x_1z_i$, $i\neq 1$. Исключительный дивизор $E$  задаётся условием $x_1=0$. Получаем следующие уравнения в этой карте  для $\tilde X$ и $E$ соответственно :
$$w^{2}= z_2+ z_3 z_4+\ldots + z_{n-1} z_{n}+\frac{1}{x_1^2}g(x_1, x_1 z_2,\ldots, x_1 z_n)$$ 
и
$$z_2=-( z_3 z_4+\ldots + z_{n-1} z_{n})+w^{2},$$ 
Таким образом,  многообразие $E$ является гладким  и рациональным (так как $E$ изоморфно аффинному пространству с координатами $z_3,\ldots z_n, w$).  Многообразие $\tilde X$ является гладким в каждой точке $E$.
\end{enumerate}
Таким образом, исключительный дивизор $E$ является гладким рациональным многообразием, следовательно,  универсально $CH_0$-тривиальным (см. \cite{CTP15} 1.5).

\qed}

\lem\label{mortr1}{Пусть $k$ -- алгебраически замкнутое поле характеристики $2$ и пусть $$X: \; y^2=f(x_1,\ldots, x_n)$$
аффинное  двойное накрытие, сингулярное в точке $P=(y, x_1,\ldots, x_n)=(0, 0,\ldots, 0)$, где $n\geq 3$ нечётно и $(0,\ldots, 0)$ -- почти невырожденная критическая точка многочлена $f$.  Пусть $\tilde X\to X$  -- раздутие точки $P$.  
Тогда :
\begin{itemize}
\item[(i)]   в окрестности точки $P$ накрытие $X$ задаётся условием 
$$y^2=ax_1^2+x_2x_3+x_4x_5+\ldots + x_{n-1}x_n+g(x_1,\ldots, x_n),$$ 
где  каждый одночлен в записи $g(x_1,\ldots,  x_n)$ имеет степень не менее трёх, и  коэффициент $b$ многочлена $g$ при $x_1^3$ не равен нулю. 
\item[(ii)]  многообразие $\tilde X$ является гладким в окрестности исключительного дивизора $E$  и многообразие $E$ является универсально $CH_0$-тривиальным.
\end{itemize}  }
\proof{Свойство $(i)$ рассматривается в упражнении V.5.7 книги \cite{Ko96} (см. также доказательство теоремы \ref{revcp} далее).  
Докажем $(ii)$.  Достаточно рассмотреть следующие карты раздутия:
\begin{enumerate}
\item  $x_i=y z_i$, $1\leq i \leq n$ и $\tilde X$ задаётся условием
$$1= az_1^2+z_2 z_3+ z_4 z_5+\ldots +  z_{n-1} z_{n}+\frac{1}{y^2}g(yz_1,\ldots, y z_n)$$ в аффинных координатах $y, z_1, \ldots z_n$. 

Заметим, что многочлен $\frac{1}{y^2}g(y z_1,\ldots, y z_n)$ делится на $y$. 
Исключительный дивизор $E$ раздутия $\tilde X\to X$ задаётся условием $y=0$. Получаем  уравнение $E$ в этой карте :
$$ az_1^2+z_2 z_3+ z_4 z_5+\ldots +  z_{n-1} z_{n}=1.$$
При $a=0$ получаем произведение $\mathbb{A}^1$ и гладкой квадрики. При $a\neq 0$ получаем неприводимую квадрику, сингулярную  в одной  точке
$z_{i}=0, i>1$ и  $az_{1}^2=1.$

Многообразие $\tilde X$ является гладким в каждой точке $E$ : сингулярная точка $\tilde X$ должна удовлетворять условиям : $z_2=\ldots= z_n=0$, $y=0$, $bz_1^3=0$ и $az_1^2=1$, что невозможно.

\item $y=wx_2$,  $x_i=x_2z_i$, $i\neq 2$. Исключительный дивизор $E$  задаётся условием $x_2=0$. Получаем следующие уравнения в этой карте  для $\tilde X$ и $E$ соответственно :
$$w^{2}= az_1^2+z_3+ z_4 z_5+\ldots + z_{n-1} z_{n}+\frac{1}{x_2^2}g(x_2z_1, x_2, x_2z_3\ldots, x_2 z_n)$$ 
и
$$z_3=-(az_1^2+ z_4 z_5+\ldots + z_{n-1} z_{n})+w^{2}.$$ 
Как и в предыдущих вычислениях,  многочлен $\frac{1}{x_2^2}g(x_2z_1, x_2, x_2z_3\ldots, x_2 z_n)$ делится на $x_2$.
Таким образом,  многообразие $E$ является гладким и рациональным (так как изоморфно аффинному пространству).   Многообразие $\tilde X$ является гладким в каждой точке $E$.

\item $y=wx_1$,  $x_i=x_1z_i$, $i\neq 1$. Исключительный дивизор $E$  задаётся условием $x_1=0$. Получаем следующие уравнения в этой карте  для $\tilde X$ и $E$ соответственно :
$$w^{2}= a+z_2z_3+ z_4 z_5+\ldots + z_{n-1} z_{n}+\frac{1}{x_1^2}g(x_1, x_1 z_2,\ldots, x_1 z_n)$$ 
и
$$a+z_2z_3+ z_4 z_5+\ldots + z_{n-1} z_{n}-w^{2}=0.$$ 
Многообразие $E$  является  неприводимой квадрикой, сингулярной в одной точке  $z_{i}=0 ;\forall i$ et $a-w^2=0$.
Так как коэффициент $g$ при $x_1^3$ не равен нулю,  многообразие $\tilde X$ является гладким в каждой точке $E$ (аналогично вычислениям в пункте $1$).

\end{enumerate}
Получаем, что многообразие $E$ является непроводимым, имеет только изолированнyю сингулярную точку $(c:1:0:\ldots 0)$, где $c^2=a$, и открытое подмногообразие $E$ является гладким и рациональным.
По лемме \ref{singord}, $E$ универсально $CH_0$-тривиально.

\qed}

\lem\label{mortr2}{Пусть $k$ -- алгебраически замкнутое поле характеристики $p>2$ и пусть $$X: \; y^p=f(x_1,\ldots, x_n)$$
аффинное  циклическое накрытие, сингулярное в точке $P=(y, x_1,\ldots, x_n)=(0, 0,\ldots, 0)$, где $(0,\ldots, 0)$ --  невырожденная критическая точка многочлена $f$.  Предполoжим, что $n$ чётно.
Тогда :
\begin{itemize}
\item[(i)]   в окрестности точки $P$ накрытие $X$ задаётся условием \\
$y^p=x_1x_2+x_3x_4+\ldots + x_{n-1}x_n+g(x_1,\ldots, x_n),$ где каждый одночлен в записи $g(x_1,\ldots,  x_n)$ имеет степень не менее трёх.
\item[(ii)]  Многообразие $\tilde X$, полученное в результате раздутия точки $P$   и  конечного числа раздутий изолированных сингулярных точек 
 с образом $P$ в $X$, является гладким в окрестности $\tilde X_P$  и слой  $\tilde X_P$  является универсально $CH_0$-тривиальным (но, в общем случае,  $\tilde X_P$  не является неприводимым).
\end{itemize}  }
\proof{Свойство $(i)$ рассматривается в упражнении V.5.6.6 книги \cite{Ko96} (см. также доказательство теоремы \ref{revcp} далее).
 Докажем $(ii)$.  
 Пусть $X'\to X$ -- раздутие $X$ в точке $P$. Достаточно рассмотреть следующие карты :
\begin{enumerate}
\item $y=wx_1$,  $x_i=x_1z_i$, $i\neq 1$. Исключительный дивизор $E$  задаётся условием $x_1=0$. Получаем следующие уравнения  в этой карте для $X'$ и $E$ соответственно :
$$x_1^{p-2}w^{p}= z_2+ z_3 z_4+\ldots + z_{n-1} z_{n}+\frac{1}{x_1^2}g(x_1, x_1 z_2,\ldots, x_1 z_n)$$
(где многочлен $\frac{1}{x_1^2}g(x_1, x_1 z_2,\ldots, x_1 z_n)$ делится на $x_1$) 
и
$$z_2=-( z_3 z_4+\ldots + z_{n-1} z_{n}),$$ 
 Таким образом,  многообразие $E$ является гладким и рациональным.  Многообразие $ X'$ является гладким в каждой точке $E$.

\item $x_i=y z_i$, $1\leq i \leq n$ и $X'$ задаётся условием
$$y^{p-2}=z_1z_2+ z_3z_4+\ldots + z_{n-1} z_{n}+\frac{1}{y^2}g(y z_1,\ldots, yz_n).$$ 
 
Исключительный дивизор $E$ раздутия $ X'\to X$ задаётся условием $y=0$. Получаем уравнение $E$ в этой карте :
$$ z_1 z_2+z_3 z_4+\ldots +  z_{n-1} z_{n}=0,$$ задающее  квадрику, сингулярную в точке $(z_1,\ldots, z_n)=(0,\ldots, 0)$. 

Многообразие $X'$
является сингулярным
 в единственной  точке $P'=(y, z_1,\ldots, z_n)=(0,\ldots 0)$ если $p>3$ и гладким, если $p=3$. Если $p>3$, пусть $X''\to X'$ раздутие $X'$ в точке $P'$. Аналогично, рассматриваем следующие карты :

\begin{enumerate}
\item $y=z_1w$, $z_i=t_iz_1$,  $i \neq 1$,
исключительный дивизор $E'$ задаётся условием $z_1=0$. Получаем следующие уравнения  в этой карте для $X''$ и $E'$ соответственно :
$$w^{p-2}z_1^{p-4}=t_2+t_3t_4+\ldots + t_{n-1}t_n+\frac{1}{z_1^2}h(z_1w, z_1, z_1t_2,\ldots, z_1t_n),$$ 
$$t_2+t_3t_4+\ldots + t_{n-1}t_n=0,$$ 
где мы обозначили 
$\frac{1}{y^2}g(y z_1,\ldots, yz_n)=h(y, z_1, \ldots, z_n).$ Заметим, что многочлен $h(z_1w, z_1, z_1t_2,\ldots, z_1t_n)$ делится на $z_{1}^3$.
Получаем, что $E'$ является гладким и рациональным: произведение $\mathbb A^1$ (соответствующего координате $w$), и многообразия
$t_2=-(t_3t_4+\ldots + t_{n-1}t_n),$
и что  $X''$ является гладким в каждой точке $E'$ в этой карте.

\item $z_i=y t_i$, $1\leq i \leq n$, $E'$ задаётся условием $y=0$ и $X''$ задаётся условием
$$y^{p-4}=t_1t_2+ t_3t_4+\ldots + t_{n-1} t_{n}+\frac{1}{y^4}g(y^2 t_1,\ldots, y^2t_n)$$ в этой карте.
Многочлен $\frac{1}{y^4}g(y^2 t_1,\ldots, y^2t_n)$ делится на $y$.
Исключительный дивизор $E'$ является квадрикой $$t_1t_2+ t_3t_4+\ldots + t_{n-1} t_{n}=0.$$

Аналогично, многообразие $X''$  является сингулярным в 
единственной
точке $(y, t_1,\ldots, t_n)=(0,\ldots 0)$ если $p>5$ и гладким, если $p=5$. 
Если $X''$  сингулярно, то мы повторяем предыдущую конструкцию. После конечного числа таких операций мы получим многообразие  $\tilde X\to X$,  гладкое в каждой точке над $P$ и такое, что все исключительные дивизоры являются рациональными многообразиями, гладкими или с единственной изолированной сингулярной точкой, 
   как описано выше.  
\end{enumerate}

\end{enumerate}

Из описания исключительных дивизоров и  леммы \ref{singord} получаем, что слой  $\tilde X_P$ является связным $CH_0$-универсально тривиальным многообразием. \\

\qed}

\lem\label{mortr3}{Пусть $k$ -- алгебраически замкнутое поле характеристики $p>2$ и пусть $$X: \; y^p=f(x_1,\ldots, x_n)$$
аффинное  циклическое накрытие, сингулярное в точке $P=(y, x_1,\ldots, x_n)=(0, 0,\ldots, 0)$, где $(0,\ldots, 0)$ --  невырожденная критическая точка многочлена $f$.  Предполoжим, что $n$ нечётно.
Тогда :
\begin{itemize}
\item[(i)]   в окрестности точки $P$ накрытие $X$ задается условием \\
 $y^p=x_1^2+x_2x_3+x_4x_5+\ldots + x_{n-1}x_n+g(x_1,\ldots, x_n)$,
где каждый одночлен в записи $g(x_1,\ldots,  x_n)$ имеет степень не менее трёх.
\item[(ii)]  Многообразие $\tilde X$, полученное в результате раздутия точки $P$   и  конечного числа раздутий изолированных сингулярных точек 

с образом $P$ является гладким в окрестности $\tilde X_P$  и слой  $\tilde X_P$  является универсально $CH_0$-тривиальным (но, в общем случае,  $\tilde X_P$  не является неприводимым).
\end{itemize}  }
\proof{Свойство $(i)$ рассматривается в упражнении V.5.6.6 книги \cite{Ko96} (см. также доказательство теоремы \ref{revcp} далее).

Докажем $(ii)$.  
 Пусть $X'\to X$ -- раздутие $X$ в точке $P$. Достаточно рассмотреть следующие карты :
\begin{enumerate}
\item $y=wx_1$,  
$x_i=x_1z_i$, $i\neq 1$.
 Исключительный дивизор $E$  задаётся условием $x_1=0$. Получаем следующие уравнения в этой карте  для $X'$ и $E$ соответственно :
$$x_1^{p-2}w^{p}= 1+z_2z_3+ z_4 z_5+\ldots + z_{n-1} z_{n}+\frac{1}{x_1^2}g(x_1, x_1 z_2,\ldots, x_1 z_n)$$
(где многочлен $\frac{1}{x_1^2}g(x_1, x_1 z_2,\ldots, x_1 z_n)$ делится на $x_1$) 
и
$$1+z_2 z_3+ z_4z_5+\ldots + z_{n-1} z_{n}=0.$$ 
 Таким образом,  многообразие $E$ является гладким и рациональным : произведение $\mathbb A^1$ (соответствующего координате $w$) и гладкой квадрики $1+z_2 z_3+ z_4z_5+\ldots + z_{n-1} z_{n}=0$.  Многообразие $ X'$ является гладким в каждой точке $E$  в этой карте.

\item $y=wx_2$,  
$x_i=x_2z_i$, $i\neq 2$.
Исключительный дивизор $E$  задаётся условием $x_2=0$. Получаем следующие уравнения   в этой карте для $X'$ и $E$ соответственно :
$$x_2^{p-2}w^{p}= z_1^2+z_3+ z_4 z_5+\ldots + z_{n-1} z_{n}+\frac{1}{x_2^2}g(z_1x_2,  x_2, x_2 z_3,\ldots, x_2 z_n)$$
(где многочлен $\frac{1}{x_2^2}g(z_1x_2,  x_2, x_2 z_3,\ldots, x_2 z_n)$ делится на $x_2$) 
и
$$z_3=-(z_1^2+ z_4 z_5+\ldots + z_{n-1} z_{n}),$$ 
 Таким образом,  многообразие $E$ является гладким и рациональным (так как изоморфно аффинному пространству).  Многообразие $ X'$ является гладким в каждой точке $E$  в этой карте.

\item $x_i=y z_i$, $1\leq i \leq n$ и $X'$ задаётся условием
$$y^{p-2}=z_1^2+z_2z_3+ z_4z_5+\ldots + z_{n-1} z_{n}+\frac{1}{y^2}g(y z_1,\ldots, yz_n).$$ 
Заметим, что многочлен $\frac{1}{y^2}g(y z_1,\ldots, yz_n)$ делится на $y$.

Исключительный дивизор $E$ раздутия $ X'\to X$ задаётся условием $y=0$. Получаем уравнения $E$ в этой карте :
$$ z_1^2+z_2z_3+ z_4z_5+\ldots +  z_{n-1} z_{n}=0,$$ задающее  квадрику, сингулярную в точке $(z_1,\ldots, z_n)=(0,\ldots, 0)$. 

Многообразие $X'$ также является сингулярным в
единственной 
 точке $P'=(y, z_1,\ldots, z_n)=(0,\ldots, 0)$ если $p>3$ и гладким в окрестности исключительного дивизора, если $p=3$. 
Если $p>3$, пусть $X''\to X'$ раздутие $X'$ в точке $P'$. Аналогично предыдущей лемме, рассматриваем следующие карты :

\begin{enumerate}
\item $y=z_1w$, $z_i=t_iz_1$, 
$i\neq 1$,
исключительный дивизор $E'$ задаётся условием $z_1=0$. Получаем следующие уравнения  для $X''$ и $E'$ соответственно :
$$w^{p-2}z_1^{p-4}=1+t_2t_3+ t_4t_5+\ldots + t_{n-1}t_n+\frac{1}{z_1^2}h(z_1w, z_1, z_1t_2,\ldots, z_1t_n),$$ 
$$1+t_2t_3+ t_4t_5+\ldots + t_{n-1}t_n=0,$$ 
где многочлен
$\frac{1}{y^2}g(y z_1,\ldots, yz_n)=h(y, z_1, \ldots, z_n)$ делится на
 $z_{1}^3$.
Получаем, что $E'$ является гладким и рациональным: произведение $\mathbb A^1$ (соответствующего координате $w$) и многообразия, заданного уравнением
 $1+t_2t_3+ t_4t_5+\ldots + t_{n-1}t_n=0.$ Многообразие  $X'$ является гладким в каждой точке $E'$  в этой карте.

\item $y=z_2w$, $z_i=t_iz_2$,
$i\neq 2$,
исключительный дивизор $E'$ задаётся условием $z_2=0$. Получаем следующие уравнения  для $X''$ и $E'$ соответственно :
$$w^{p-2}z_2^{p-4}=t_1^2+t_3+ t_4t_5+\ldots + t_{n-1}t_n+\frac{1}{z_2^2}h(z_2w, z_2t_1, z_2, z_2t_3,\ldots, z_2t_n),$$ 
$$t_1^2+t_3+ t_4t_5+\ldots + t_{n-1}t_n=0,$$ 
где мы обозначили 
$\frac{1}{y^2}g(y z_1,\ldots, yz_n)=h(y, z_1, \ldots, z_n).$
Получаем, что $E'$ является гладким и рациональным и что $X''$  является гладким в каждой точке $E'$ в этой карте.

\item $z_i=y t_i$, $1\leq i \leq n$, $E'$ задаётся условием $y=0$ и $X''$ задаётся условием
$$y^{p-4}=t_1^2+t_2t_3+ t_4t_5+\ldots + t_{n-1} t_{n}+\frac{1}{y^4}g(y^2 t_1,\ldots, y^2t_n).$$
Многочлен $\frac{1}{y^4}g(y^2 t_1,\ldots, y^2t_n)$ делится на $y$.
Исключительный дивизор $E'$ является квадрикой $$t_1^2+t_2t_3+ t_4t_5+\ldots + t_{n-1} t_{n}=0.$$

Аналогично, многообразие $X''$  является сингулярным в окрестности $E$ в этой карте в единственной точке
$(y, t_1,\ldots, t_n)=(0,\ldots 0)$ если $p>5$ и гладким, если $p=5$. 
Если $X''$ сингулярно, то мы повторяем предыдущую конструкцию. После конечного числа таких операций мы получим многообразие  $\tilde X\to X$,  гладкое в каждой точке над $P$ и такое, что все исключительные дивизоры являются рациональными многообразиями, гладкими или с единственной изолированной сингулярной точкой, 
как описано выше. 
\end{enumerate}

\end{enumerate}

Из описания исключительных дивизоров и леммы \ref{singord} получаем, что слой  $\tilde X_P$ является связным $CH_0$-универсально тривиальным многообразием. \\

\qed}

\vspace{1cm}

Следующее утверждение даёт ключевые нетривиальные инварианты циклических накрытий.

Напомним, что коэффициенты многочленов  $f\in k[x_0,\ldots, x_n]$ заданной степени  параметризуются точками некоторого аффинного  пространства.  Общий выбор коэффициентов $f$ означает,  что мы рассматриваем коэффициенты из некоторого непустого открытого (в топологии Зарисского) подмногообразия этого аффинного пространства.\\

\theo\label{revcp}{
Пусть $k$ -- алгебраически замкнутое поле характеристики $p$ и пусть $f(x_0,\ldots, x_n)$  -- однородный многочлен степени $mp\geq n+1$, $n\geq 3$. Для общего 
выбора коэффициентов $f$ выполняются следующие условия:
\begin{itemize}
\item[(i)] все критические точки $f$ являются невырожденными, если $p>2$  или $p=2$ и $n$ чётно;
\item[(ii)]  все критические точки $f$ являются  почти невырожденными, если  $p=2$ и $n$ нечётно.
\item[(iii)]  Если $
\tilde X\to X$ -- разрешение особенностей $X$, полученное последовательным раздутием сингулярных точек, то морфизм $\tilde X\to X$ является $CH_0$-универсально тривиальным,
  $H^0(\tilde X, \Lambda^{n-1}\Omega_{\tilde X})\neq 0$ и многообразие $\tilde X$ {\it не является} $CH_0$-универсально тривиальным.  \end{itemize}
}

\proof{ Свойства (i) и (ii) следуют из упражнений \cite{Ko96}V.5.7 и 5.11.  Предположим, $P=(b, a_1,\ldots, a_n)$ -- критическая точка $f$. При помощи линейной замены переменных $y-c, x_i-a_i$, где $c^p=f(P)$ (поле $k$ алгебраически замкнуто)  можно предположить, что $P=(0,\ldots 0)$.
 Тогда можно разложить $f$ в сумму линейной части, квадратичной части и части высших степеней: $f=f_1(x_1,\ldots x_n)+f_2(x_1,\ldots x_n)+f_3(x_1,\ldots x_n).$ Так как $P$ -- критическая точка, то $f_1=0$. Так как поле $k$ алгебраически замкнуто, то каждая квадратичная форма над $k$ может быть представлена в диагональной форме либо суммой квадратов  (если характеристика поля не  равна $2$), либо суммой $\sum x_iy_i$ (регулярная часть) и суммой квадратов.  
 Таким образом, несложно проверить,
  что условие, что $P$ невырожденная (соотв. почти невырожденная) точка, соответствует разложению в леммах 
 \ref{mortr}, \ref{mortr2}, \ref{mortr3} (соотв. \ref{mortr1}), что выполняется для общего выбора коэффициентов $f$ (см. упражнение \cite{Ko96} V. (5.6.6.3)).  

Для доказательства $(iii)$,  как и в работе Тотаро \cite{To15}, мы используем  \cite{Ko96} V.5.11 для $ \mathbb P^n_k,$ $n\geq 3$ и $ L^p=\mathcal O_{\mathbb P^n}(mp)$. Получаем :

\begin{enumerate}
\item  $K_{\mathbb P^n}\otimes L^p=\mathcal O_{\mathbb P^n}(mp-n-1)$,
\item   При условии $mp\geq  4$ 
отображение $H^0({\mathbb P^n},L^p)\to \mathcal O_{{\mathbb P^n}}/m_x^4\otimes L^p$ сюрьективно для любой замкнутой точки $x\in X$.
\end{enumerate}

Kaк следует из \cite{Ko96} V.5.7 (см. также  \cite{Ko96} V.5.11,  \cite{KSC04} Теорема 4.4),
для общего выбора $f\in H^0(\mathbb P^n_k, \mathcal O_{\mathbb P^n}(mp))$ (в частности, удовлетворяющего $(i)$ и $(ii)$), если $q: X\to \mathbb P^n_k$ -- циклическое накрытие $\mathbb P^n_k$ степени $p$, разветвлённое над гиперповехностью $f=0$ и  $\pi: \tilde X\to X$ -- разрешение особенностей $X$, полученное последовательным раздутием сингулярных точек, то $\pi^*q^*\mathcal O_{\mathbb P^n}(mp-n-1)$ является подпучком  $\Lambda^{n-1}\Omega_{\tilde X}$,  
в частности, 
 если   $mp-n-1\geq 0$,  то
$H^0(\tilde X, \Lambda^{n-1}\Omega_{\tilde X})\neq 0$. 

  Как доказано в работе Тотаро \cite{To15} (Лемма 2.2),   если $\tilde X$ является $CH_0$-универсально тривиальным, то $H^0(\tilde X, \Lambda^{n-1}\Omega_{\tilde X})=0$. 
 Из лемм \ref{mortr}, \ref{mortr1}, \ref{mortr2}, \ref{mortr3}  следует, что все слои морфизма $\tilde X\to X$  являются $CH_0$-универсально тривиальными, значит, и морфизм  $\tilde X\to X$  также $CH_0$-универсально тривиальный (см. \cite{CTP15} 1.8)
 \qed\\ 

 {\it Замечание}. Если $n+1>mp-m$  и многообразие $X$ нормально (что верно, в частности, если $f$ имеет только  изолированные критические точки), то $X$ является многообразием Фано : пучок $-K_X$ обилен (см. \cite{KSC04} 4.14). \\

 \section{Накрытия, которые не являются стабильно рациональными}
 
 \theo\label{revnsr}{Пусть $p$ -- простое число. Пусть $X$ -- циклическое накрытие $\mathbb P_{\mathbb C}^n, n\geq 3$ степени $p$,   разветвлённое над очень общей гиперповерхностью $f(x_0,\ldots , x_n)=0$ степени $mp$.   Предположим,  что 
 $m(p-1) <n+1\leq mp$.  Тогда $X$  -- многообразие Фано, которое не является стабильно рациональным  многообразием.  Существует не стабильно рациональное накрытие   степени $p$,   разветвлённое  над гиперповерхностью  степени $mp$, определённое над числовым полем.  
 \proof{  
 Пусть $Y:  y^p=f(x_0,\ldots ,x_n)$ -- накрытие, удовлетворяющее условиям теоремы  \ref{revcp}.  Можно выбрать $Y$ таким образом, что коэффициенты многочлена $f$ определены над некоторым конечным  полем $\mathbb F_q$ характеристики $p$;  так как в  теореме  \ref{revcp}   условие $f=0$ задаёт очень общую гиперповерхность, то  можно также предположить, что гиперповерхность  $f=0$ является гладкой над   $\mathbb F_q.$ Таким образом,  существует многочлен  $H$ степени $mp$ с коэффициентами в некотором числовом поле, такой, 
   что $f$ является редукцией многочлена $H$ по модулю $p$,  и что   
    накрытие  $X : y^p=H(x_0, \ldots, x_n)$ 
является гладким многообразиeм. 
 Так как $X$ вырождается в $Y$ и разрешение особенностей $Y'\to Y$, построенное в леммах \ref{mortr}, \ref{mortr1}, \ref{mortr2},
\ref{mortr3}  и теореме \ref{revcp}  является $CH_0$-универсально тривиальным,  то $X_{\mathbb C}$ не является  $CH_0$-универсально тривиальным многообразием по теореме \ref{revcp} и  \cite{CTP15} 1.14(iii). Следовательно, $X$ не является стабильно рациональным. Более того, по построению  мы получаем пример над числовым полем. 
Так как коэффициенты многочленов степени $mp$ с комплексными коэффициентами параметризуются неприводимым многообразием, то по   \cite{CTP15} 2.3 для очень общего выбора такого многочлена, соответствующее
 накрытие $\mathbb P_{\mathbb C}^n$ степени $p$ не является $CH_0$-универсально тривиальным. \\
 
 \qed}
 
  {\it Замечание}.  Таким образом, мы получаем, что циклическое накрытие $X$ проективного пространства $\mathbb P_{\mathbb C}^n$ степени $p$,   разветвлённое над очень общей гиперповерхностью $f(x_0,\ldots , x_n)=0$ степени $mp$ не является $CH_0$-универсально тривиальным.  Аналогично работе \cite{CTP15}, из этого следует, что $X$  не является рациональным ретрактом. Напомним, что стабильно рациональное многообразие является рациональным ретрактом, однако до сих пор неизвестно, различны ли эти понятия.  \\
 
 \vspace{1cm}
 
 {\bf Примеры.}
 \begin{enumerate}
 \item При $p=2, n=3$ и $mp=6$ получаем другое доказательство результата А. Бовиля \cite{B6}.
 
 \item При $n=3, m=p=2$  получаем  двойные  накрытия $\mathbb P^3_{\mathbb C}$, разветвлённые над квартикой (более общие результаты получены в работе  Клэр Вуазан \cite{V13}).

\item При $p=2, n=4$ и   $mp=6, 8$ получаем, что двойное накрытие 
$\mathbb P^4_{\mathbb C}$, разветвлённое над очень общей гиперповерхностью степени $6$ или $8$ не является стабильно рациональным.
\item При $p=2, n=5$ получаем примеры для  $2m= 8, 10$.

\item При $p=3, n=4$ и $mp=6$ получаем пример многообразия Фано размерности $4$, которое не является стабильно рациональным : накрытие степени $3$, разветвлённое над очень общей поверхностью степени $6$.

 \item
 Результаты о двойных накрытиях
 $\mathbb P^n_{\mathbb C}$,  $n=4,5$,  разветвлённых над квартикой
 (А. Бовиль \cite{B4}), не покрываются результатами теоремы \ref{revnsr}.

 \end{enumerate}

\end{document}